\documentstyle{article}

\def\ben{\begin{equation}}
\def\een{\end{equation}}
\def\bea{\begin{eqnarray}}
\def\eea{\end{eqnarray}}
\begin{document}

\title{Gravitational Instantons, Confocal Quadrics
and Separability of the Schr\"odinger and Hamilton-Jacobi equations.}
\author{G. W. Gibbons
\\
C.M.S,
\\ Cambridge University,
\\ Wilberforce Rd.,
\\ Cambridge CB3 OWA,
 \\ U.K.}
\maketitle

\begin{abstract}A hyperk\"ahler 4-metric  with
a triholomorphic $SU(2)$ action gives rise to
a family of confocal quadrics in Euclidean 3-space
when cast in
the canonical form of a hyperk\"ahler 4-metric
metric with a triholomorphic circle action.
Moreover, at least in the case of geodesics orthogonal to
the $U(1)$ fibres, both the covariant Schr\"odinger and
the Hamilton-Jacobi equation is
separable and the system integrable.

\end{abstract}

\section{ Introduction}
The work described in what follows was done some time ago, and arose
as an attempt to understand in a more conceptual way a calculation
given in the appendix of \cite{GORV}. The reason for returning to it
here is two-fold.
\medskip \begin{itemize} \item Its relevance to some recent recent work of
Dunajski \cite{Dunajski} on quadrics in general and integrability.
The example discussed below is particularly  striking
in that respect because one may relate it to a problem
of genuine practical significance.

 \item It's relevance to the problem of consistent
sphere reductions in Kaluza-Klein supergravity and
string theories in various dimensions for which the properties
of quadrics are essential \cite{CGHP}. The present work is concerned
exclusively with particular  four-dimensional riemannian metrics which
can be seen to exhibit some of the general features discussed in
\cite{CGHP} in a particularly intriguing way.

\end{itemize}

\section {Left and Right invariant vector
fields}

What follows is a brief summary of some properties
of $SU(2)$ which will be needed in the sequel.
Suppose that $\sigma^i$ are left-invariant one forms with
\ben
d \sigma ^1 = \sigma ^2 \wedge \sigma ^3 \quad {\it etc.}
\een
Then the dual basis of left-invariant vector fields ${\bf L}_i$
satisfies
\ben
\bigl [ {\bf L}_1, {\bf L}_2 \bigr ] =-{\bf L}_3 \quad {\it etc.}
\een
The generators of left translations form a right-invariant basis ${\bf R}_i
$
which commutes with the left-invariant basis:
\ben
\bigl [ {\bf L}_i, {\bf R}_j \bigr ]=0 \label{commute},
\een
and satisfies:
\ben
\bigl [ {\bf R}_1, {\bf R}_2 \bigr ] ={\bf R}_3 \quad {\it etc.}
\een
Now the left and right invariant vector fields are not linearly independent
and hence
\ben
\pmatrix{ {\bf R}_1 \cr {\bf R}_2 \cr {\bf R}_3 } =
\pmatrix {l_1 &l_2 & l_3 \cr m_1 & m_2 & m_3 \cr n_1 & n_2 & n_3 \cr}
\pmatrix{ {\bf L}_1 \cr {\bf L}_2 \cr {\bf L}_3 },
\een
where
\ben
\pmatrix {l_1 &l_2 & l_3 \cr m_1 &m_2 & m_3 \cr n_1 & n_2 & n_3 \cr}
\een
is an orthogonal matrix which depends upon the Euler angles of
$SU(2)$. Thus for example:
\ben
l_1^2 + l_2 ^2 + l_3 ^2 =1. \label{norm}
\een
Moreover (\ref{commute}) implies that
\ben
{\bf L}_i l_j =- \epsilon _{ijk} l_k \label{diff},\een
which may also be written as
\ben
dl_i=- \epsilon _{ijk} l_j \sigma ^
k \label{diffy}.
\een

\section{ BGPP  metrics}

If $SU(2)$ acts triholomorphically  the metric may be cast in the form:
\ben
ds^2 = ABC d \eta ^2 + {BC \over A}( \sigma ^1) ^2 +{CA  \over B}( \sigma
 ^2) ^2+{AB \over C}( \sigma ^3) ^2
\een
where $\eta$ labels the orbits of $SU(2)$. If dot
indicates differentiation with respect to $\eta$,
then  $A,B,C$ satisfy
\ben
{\dot A} = BC \quad {\it etc.} \label{eqns}
\een
The  K\"ahler forms $\Omega_i$
are given by
\ben
\Omega_1 = BC d \eta \wedge \sigma ^1 + A \sigma ^2 \wedge \sigma ^3 \quad
 {\it etc.}
\een
Note that (\ref{eqns}) is equivalent to the closure of the
$\Omega_i$.
For future reference, we record that
\ben
\Omega_1 ({\bf R}_1,)=-BCl_1 d \eta + A l_2 \sigma ^3 - A l_3 \sigma ^2 \label{inner}.
\een

Because $\langle {  \sigma ^j, \bf L}_i  \rangle = \delta ^j_i$,
one easily sees that
\ben
g( {\bf L}_1, {\bf L}_1 )= { BC \over A} \quad {\it etc} \label{metric}.
\een

The equations (\ref{eqns}) are of a similar  form to Euler's equations for
a rigid body but they do not correspond to any set of moments of inertia.
Nevertheless the equations may still be   integrated  completely.
In order to avoid elliptic functions we
change the variable labelling the orbits of $SU(2)$.
Thus we define  $\lambda$ by
\ben
ABC d \eta ^2 = { 1\over 4} { d \lambda ^2 \over ABC}.
\een
and find that
\ben
A= \sqrt{\lambda - \lambda _1} \quad {\it etc.}  \label{soln}
\een

\section {Canonical form}

The canonical form for a self-dual metric with a triholomorphic $U(1)$ action is (see e.g. \cite{GR2})
\ben
ds^2  = V^{-1} (d \tau + \omega_i dx^i)^2 + V ( dx^2 +dy^2 + dz ^2)
\label{Harm},
\een
where
\ben
{\rm curl} \quad \omega = {\rm grad} \quad  V\label{curl}.
\een

It follows that the  function  $V$ is
harmonic on the Euclidean space  whose Cartesian coordinates
are the three moment maps $x,y,z$. This Euclidean space is
of course just the quotient of the 4-manifold
by the triholomorphic  circle action generated by $\partial / \partial \tau
$.
Moreover it is clear that given $V$ one may
reconstruct the metric.

The three  K\"ahler forms are
\ben
\Omega _1= (d \tau + \omega_i dx^i) \wedge dx - V d y \wedge dz \quad{\it etc} \label{forms}
\een
Note that the closure of the $\Omega_i$ follows from (\ref{curl}).

From (\ref{forms}) it follows that $(x,y,z)$ are the three moment maps
for the Killing field $\partial /\partial \tau$, i.e.
\ben
\Omega_1( {\partial \over \partial \tau} , ) = d x \quad {\it etc}.
\een
Moreover we have the obvious relation:
\ben
V^{-1} = g({\partial \over \partial \tau},{\partial \over \partial \tau}).
\een

\section{Comparison}

We pick one of the generators, ${\bf R}_1$ say, of left translations.
This generates a $U(1)$ subgroup of $SU(2)$
which we identify with $\partial/\partial \tau$.
Now if
\ben
x=-l_1 A \quad {\it etc},
\een
then by (\ref{eqns}) and (\ref{diffy})
\ben
dx= -l_1 BC d \eta +A l_2 \sigma ^3 - A l_3 \sigma ^2.
\een
but using (\ref{inner}) we have
\ben
dx=\Omega _1({\bf R}_1,\thinspace ).
\een
Thus we have identified the moment maps.
Using (\ref{norm}) and (\ref{soln}) we deduce that
\ben
{ x^2 \over \lambda - \lambda _1} +{ y^2 \over \lambda - \lambda _2}+
{ z^2 \over \lambda - \lambda _3}=1.
\een
The cosets  $SU(2)/U(1)$ appear as set of confocal
quadrics in Euclidean 3-space, labelled by the variable $\lambda$.
Each non-degenerate orbit is a circle bundle over the
quadric. As  bundles they are precisely the Hopf bundle.
It remains to find the harmonic function $V$. This is of course
given by
\ben
V^{-1} = g({\bf R}_1, {\bf R}_1).
\een
Because
\ben
{\bf R}_1= l_1 {\bf L}_1 +l_2 {\bf L}_2+l_3 {\bf L}_3,
\een
one has
\ben
V^{-1} = \sqrt{(\lambda -\lambda _1)(\lambda -\lambda _2)(\lambda -\lambda _3)} \Bigl ( {x^2 \over (\lambda - \lambda _1)^2} + {y^2 \over (\lambda -
 \lambda _2)^2}+{z^2 \over (\lambda - \lambda _3)^2}  \Bigr ).
\een

The result just described was obtained originally in a rather less
conceptual way
 in \cite{GORV}. Another approach is described in \cite{Dunajski}.

It is  a straightforward task to express this expression
in ellipsoidal coordinates and check explicitly that it is indeed
a solution of Laplace's equation. By
doing so we obtain a bonus: we discover that, in a special
case at least, we
 may separate the Hamilton-Jacobi equation governing the geodesic flow.

\section{ Ellipsoidal Coordinates}

We begin by introducing ellipsoidal coordinates in
three dimensional Euclidean space in
the usual way.
We assume that  $\lambda _1 < \lambda _2 < \lambda _3$
and  restrict $\lambda$ to the range $ \lambda _3 \le \lambda \le \infty$.
Two further coordinates $\mu$ and $\nu$ lie in the intervals $\lambda \le \mu\le \lambda_2 $ and $\lambda_2 \le \nu \le \lambda _3$.
The Cartesian coordinates are given by
\ben
x^2= { ( \lambda -\lambda _1) ( \mu-\lambda_1) ( \nu- \lambda_1) \over
(\lambda _1 - \lambda _2 ) ( \lambda _1-\lambda _3)},
\een
\ben
y^2={ ( \lambda -\lambda _2) ( \mu-\lambda_2) ( \nu- \lambda_2) \over
(\lambda _2 - \lambda _1 ) ( \lambda _2-\lambda _3)},
\een
\ben
z^2={ ( \lambda -\lambda _3) ( \mu-\lambda_3) ( \nu- \lambda_3) \over
(\lambda _1 - \lambda _2 ) ( \lambda _2-\lambda _3)}.
\een

If we define
\ben
R(\lambda)= \sqrt {(\lambda -\lambda_1)(\lambda -\lambda_2)(\lambda -\lambda_3)},
\een
we find, after a short calculation, that
\ben
V= {  R(\lambda) \over (\lambda -\mu) ( \lambda -\nu)}.
\een

By writing out the flat space Laplace operator in
ellipsoidal coordinates, one checks that indeed $V$ is harmonic.
In fact $V$ coincides with the harmonic
function $S_\lambda$ defined in \cite{SM}.

\section{The one-form $\omega_i d x^i$}

The harmonic property of the function $V$ guarantees
that a local solution of (\ref{curl}) exists, but finding it may
not be easy. Rather remarkably,
an explicit solution may be found.

If we define
\ben
S(\mu) = \sqrt{-(\mu -\lambda _1)(\mu - \lambda _2)( \mu- \lambda _3) },
\een
and
\ben
T(\nu) = \sqrt{(\nu -\lambda _1)(\nu - \lambda _2)( \nu- \lambda _3) },
\een
then, according to  Harry Braden
one may choose:
\ben
\omega_i d x^i={ 1\over 2}
 { (\lambda-\nu) \over (\mu-\lambda) (\mu- \nu)} { S(\mu) \over T(\nu)} d \nu + { 1\over 2} { (\lambda-\mu) \over (\nu-\lambda) (\nu -\mu) } { T(\nu) \over S(\mu)} d \mu.
\een

\section{Special cases}

If all the constants $\lambda_i$ are equal  we obtain flat Euclidean
4-space.  In this case, up to scaling,
\ben
V = { 1\over r},
\een
where $r$ is the distance to the origin.
If two of the constants
coincide, one obtains the Eguchi-Hanson metric.
This is a complete metric on the cotangent bundle
of the two-sphere, $T^\star(CP^1)$, and is in fact invariant
under the action of $U(2)$.  In this case, up to scaling,
\ben
V = { 1\over r_1 } + { 1 \over r_2 },
\een
where $r_1$ and $r_2$ are the distances form two fixed points.

If the three
constants are distinct, we obtain a metric with singularities.
A simple description of the potential $V$ would be
interesting. Its singularities lie on the focal conics.
Another interesting question is what, if anything,
is the relation of these potentials
to those described in appendix 15 of Arnold's
textbook on dynamics.

\section {Integrability
of the geodesic flow}

The geodesic flow on the tangent bundle of the
general metric with triholomorphic $U(1)$ can be subjected
to a Hamiltonian reduction to get
a mechanical system on Euclidean space. This is most
simply done using the Hamilton-Jacobi equation.
One may separate off the $\tau$  dependence by writing
the action as
\ben
S= e \tau + W(x^i).
\een
One gets :
\ben
e^2 V+ V^{-1} (\nabla _i W- e \omega_i) ^2=1.
\een
If the conserved momentum $e$ vanishes, which means that the
geodesics are orthogonal to the $U(1)$ fibres, then
we obtain the Hamilton-Jacobi equation for
for a particle in an  (attractive)
gravitational  potential $V$ moving
with zero total energy. If the momentum
$e$ is not zero, then  there is an additional contribution
to the potential and also a velocity dependent
magnetic force
given by $\omega_i$.

It has been known since the time of Euler that
the motion of a planet  moving around
two  fixed gravitating centres
is completely integrable. It turns out that the
geodesic flow  in the Eguchi  metric is
completely integrable. Thus Euler's result may be
generalized to this case. Moreover
the Schr\"odinger  equation also separates \cite{GR2,M,Malmedier}.
A rather full discussion is
given in \cite{M, Malmedier}. Shortly we shall see that
at least for a special class of geodesics,
the Hamilton-Jacobi equation separates.

Before doing so we shall make some remarks.

\begin{itemize}

\item It is easy to see that integrability does not follow just from
the $SU(2)$ symmetry, one needs a further  commuting quantity.

\item Nevertheless, it has been shown that
the gravitational instanton constructed
from the helicoid using the J\"orgens correspondence
between solutions of the Monge-Amp\'ere equation and minimal surfaces
admits the tri-holomorphic action of
the Bianchi group $VII_0$, i.e. the Poincar'e
group of two dimensions $E(1,1)$.
Moreover the closely related metric, obtained from the catenoid
admits the triholomorphic action
of the Bianchi group $VI_0$, that is the Euclidean group
in two dimensions $E(2)$. In both cases the authors  find
that they can separate both the Hamilton-Jacobi equation and the
Laplace equation. Since $E(2)$ is a contraction of $SU(2)$
the catenoid metric is a  are limiting case
of the BGPP metrics. Thus it is not unreasonable
to look for separability in our case.

\item On the negative side, according to a theorem quoted in the last
section of \cite{GGR}
the planar many centre problem with more than
two fixed centres is not integrable. However,
multi-centre HyperK\"ahler  metrics are not, for
more than two centres, associated with quadrics.

\item According to \cite{FO} the general {\sl repulsive}
 Coulomb problem
motion
with {\sl freely moving}  centres is Liouville integrable.
Our forces however are {\sl attractive} and the centres are {\sl fixed}.

\end{itemize}

The issue can clearly only be settled by a calculation.

\section{Separability of the Hamilton-Jacobi Equation}

Consider now the geodesics
with vanishing momentum along the $U(1)$ fibres.
As explained above, the Hamilton-Jacobi equation becomes
\ben
|\nabla W |^2= V.
\een
In ellipsoidal coordinates this becomes
\ben
{ 1 \over h_\lambda ^2} \bigl ( {\partial W \over \partial \lambda} \bigr )
^2
+{ 1 \over h_\mu ^2} \big ( {\partial W \over \partial \mu}
\bigr )^2 +{ 1 \over h_\nu ^2} \bigl ( {\partial W \over \partial \nu}\bigr
 )^2  ={R(\lambda)
\over (\lambda -\mu) (\lambda -\nu)}.
\een

The functions $h^2_\lambda, h^2_\mu, h^2_\nu$ appear in the flat
Euclidean metric which in ellipsoidal coordinates is given by
\ben
ds^2= h^2_\lambda d \lambda ^2 + h^2_\mu d \mu ^2 + h^2 _\nu d \nu ^2.
\een

If, as before,  we define
\ben
S(\mu) = \sqrt{-(\mu -\lambda _1)(\mu - \lambda _2)( \mu- \lambda _3) },
\een
and
\ben
T(\nu) =\sqrt{(\nu -\lambda _1)(\nu - \lambda _2)( \nu- \lambda _3) },
\een
we have
\ben
2h_\lambda= { \sqrt{(\lambda -\mu)( \lambda -\nu)} \over R(\lambda)},
\een

\ben
2h_\mu= { \sqrt{(\lambda -\mu)( \mu -\nu)} \over S(\mu)},
\een
and

\ben
2h_\nu={ \sqrt{(\nu -\mu)( \nu-\lambda ) }\over T(\nu)}.
\een
It is now a simple task to check that the Hamilton-Jacobi equation
with vanishing charge $e$ separates. If however $e\ne0$ the term
involving the vector potential enters the equation and so far
I have had no success in separating it.

\section {Integrability of the Schr\"odinger Equation}

In this section we shall show that
not only the Hamilton-Jacobi equation but
also the Schr\"odinger equation is admits separable solutions
in the case of vanishing charge $e$.
In fact, by the WKB approximation, this makes the former
property almost obvious. In fact since,
by virtue of the self-duality property,
the metric
admits two covariantly constant spinor fields, the eigen functions
of other operators may  be obtained from the scalar
eigen functions in the manner described by Hawking and Pope
for the metric on K3 \cite{HP}. Thus the separability
of the Laplace operator is the key to understanding
other separability properties.

The covariant Laplace or Schr\"odinger equation:

\ben
-\nabla ^\alpha \nabla _\alpha \Psi=E \Psi
\een
becomes, when written out for metric (\ref{Harm}),
\ben
V^{-1} (\nabla_k- \omega _k{\partial \over \partial \tau}) (\nabla^k- \omega ^k{\partial \over \partial \tau}) \Psi + V {\partial ^2 \over \partial \tau ^2}
\Psi = -E \Psi.
\een
We may separate off the
$\tau$ dependence by setting $\Psi=\exp( {ie\tau}) \psi$
and obtain
\ben
V^{-1} (\nabla_k- ie \omega _k) (\nabla^k
- ie \omega ^k) \psi -e^2 V  \psi =-E \psi.
\een
In the case that the charge $e$ vanishes we get
\ben
\nabla _k \nabla ^k \psi =- E V \psi \label{shr}.
\een

Separation of variables in ellipsoidal
coordinates for equation (\ref{shr})
(without assuming that the function $V$ is harmonic)
is discussed
in \cite{MF} where a sufficient condition
is given on the function $V$.

This condition is
\ben
V={ f_\lambda \over  h^2_\lambda} + {f_\mu \over h^2_\mu} + { f_\nu  \over h^2_\nu},
\een
where $(f_\lambda, f_\mu, f_\nu)$ are arbitrary functions of
$(\lambda, \mu ,\nu)$ respectively.
The sufficient condition is satisfied with
the only one of the three functions being
non-vanishing:
\ben
f_\lambda= { 1 \over 4 R(\lambda)}.
\een
in our case.

As far as I am aware, there is no discussion of the separability
of the Schr\"odinger equation with a vector potential.

\section{Acknowledgements} It is a pleasure to thank
Harry Braden, Galliano Valent and Maciej Dunajski for helpful
discussions. Some of this work was carried out during a visit to
I.H.E.S. I would like to thank the Director and his staff for their
hospitality.


\begin{thebibliography}{99}


\bibitem{Dunajski} M Dunajski, Harmonic functions, central quadrics
and twistor theory, DAMTP preprint
\bibitem{CGHP} M Cvetic, G W Gibbons H L\"u and C N Pope, unpublished
work

\bibitem{BGPP} V Belinskii, G W Gibbons, D N Page \& C N Pope,
Asymptotically Euclidean Bianchi IX Metrics in Quantum Gravity
{\sl Phys. Lett.} {\bf B76} 433-435 (1978)

\bibitem{GR1} G W Gibbons \& P J Ruback,
The Hidden Symmetries of Taub-NUT and Monopole Scattering
 {\sl Phys. Lett.} {\bf B188}  226--230 (1987)

\bibitem{GR2} G W Gibbons \& P J Ruback,
The Hidden Symmetries of Multi Centre Metrics
{\sl Comm. Math. Phys.} {\bf 115} 267--300 (1988)


\bibitem{GORV} G W Gibbons,  D Olivier, G Valent \& P J
Ruback, Multi-Centre Metrics and Harmonic Superspace
 {\sl Nucl. Phys.} {\bf B296} 679--696 (1988)

\bibitem{M} S Mignemi,
Classical and Quantum Motion in Eguchi-Hanson Space
{\sl J Math Phys} {\bf 32} 3047-3084 (1991)

\bibitem{Malmedier} A Malmedier, The eigenvalue equation on
Eguchi-Hanson space {\tt arXiv: math.DG/0210081}

\bibitem{FO} G Fusco and W M Oliva, Integrability of a
system of N electrons Subjected to Coulombian Interactions, {\sl J Diff Eqns} {\bf 135} 16-40 (1997)

\bibitem{GGR} G W Gibbons, R Goto and P Rychenkova, HyperK\"ahler
Quotient Construction of BPS Monopole moduli Spaces, {\sl Comm Math Phys }{\bf 185} 581-599 (1997)

\bibitem{SM} S J Madden,
A Separable Potential in Triaxially Ellipsoidal Coordinates satisfying the
Laplace equation,
{\sl Celestial Mechanics} {\bf 2} 21--227 (1970)

\bibitem{AHKN} A N Aliev, M Hartascu, J Kalayci \& Y Nutku
, Gravitational Instantons from Minimal Surfaces,  {\sl Class Quant Grav}
{\bf 16} (1999) 189-210 {\tt arXiv: gr-qc/9812098}
\bibitem{HP} S W Hawking \& C N Pope, Symmetry Breaking
by Instantons in Supergravity, {\sl Nucl Phys} {\bf B146} 381 (1978)

\bibitem{MF} P M Morse \& H Feshbach {\sl Methods of Theoretical Physics}
(New York, McGraw-Hill) (1953)

\end{thebibliography}
\end{document}